\documentclass{article}
\input amssym.def\input amssym
\usepackage{graphics}

\newtheorem{thm}{Theorem}
 
 \newtheorem{lem}[thm]{Lemma}

\newtheorem{alg}[thm]{Algorithm}

\begin{document}

\title{DGMRES method augmented with eigenvectors for computing the Drazin-inverse solution of  singular
linear systems  }

\author{MENG Bin \\
{\small College of Science, Nanjing University of Aeronautics and
Astronautics,} \\ {\small Nanjing 210016,  People's Republic
 of China}\\ {\small E-mail: b.meng@nuaa.edu.cn}}

 \date{}
\maketitle

\noindent{\bf Abstract.} The DGMRES method for solving
Drazin-inverse solution of  singular linear systems is generally
used with restarting. But the restarting often slows down the
convergence and DGMRES often stagnates. We show that adding some
eigenvectors to the subspace can improve the convergence just like
the method proposed by R.Morgan in [R.Morgan, A restarted GMRES
method augmented with eigenvectors, SIAM J.Matrix Anal.Appl. 16
(1995)1154-1171]. We derive the implementation of this method and
present some numerical examples to show the advantages of this
method.

\noindent {\bf Mathematics Subject Classification.} 65F15, 15A18

\noindent{\bf Key words and phrases.} Drazin-inverse; DGMRES; Krylov
subspace; Iterative method; Eigenvector.

\section{Introduction}
DGMRES is a new iterative method for computing the Drazin-inverse
solution of linear systems \cite{sid2}. Consider the linear system
\begin{equation}\label{1}
Ax=b
\end{equation}
where $A\in\Bbb{C}^{n\times n}$ is a singular matrix. We recall that
the Drazin-inverse solution of (\ref{1}) is the vector $A^Db$, where
$A^D$ is the Drazin-inverse of the singular matrix $A$.

In \cite{sid1}, A.Sidi proposed a general approach to Krylov
subspace methods for computing Drazin-inverse solution. In that
paper the authors do not put any restriction on $A$, that is, A is
non-hermitian or hermitian, index of $A$ is arbitrary and the
spectrum of A can have any shape. In \cite{sid2}, A.Sidi gave one of
the Krylov subspace method named DGMRES, which is a GMRES-like
algorithm. Like GMRES, in practical use, we often propose restarted
DGMRES which denoted by DGMRES(m). DGMRES(m) is an economical
computing and storgewise method for the Drazin-inverse solution. But
restarting slows down the convergence and often stagnates (see
\cite{zhou1},\cite{zhou2}). In the present paper we propose a DGMRES
method augmented with eigenvectors, which can accelerate the
convergence and overcome the stagnation. Classical GMRES augmented
with eigenvectors was studied by R.Morgan (see \cite{mor}). Now we
derive a DGMRES  augmented with eigenvectors similarly. We give the
convergence analysis of DGMRES augmented with eigenvectors which
shows this method is more effective than DGMRES without augment.
Then some numerical experiments are presented to show the
convergence rate of  DGMRES augmented with eigenvectors is
remarkably improved, especially when the matrix has a very large or
small nonzero  eigenvalue.

The paper is organized as follows. In section 2, we will give a
brief review of DGMRES. In section 3, we obtain the convergence
analysis of DGMRES augmented with eigenvectors and derive the
algorithm. In section 4, we present some numerical examples to
compare the DGMRES with the DGMRES augmented with eigenvectors.

\section{DGMRES}
Throughout this paper, we suppose the index of $A$ is known. The
index of $A$, denoted by $ind(A)$, is the size of the largest Jordan
block corresponding to the zero eigenvalue of $A$. DGMRES is a
Krylov subspace methods for computing the Drazin-inverse solution
$A^Db$. For more details we refer the readers to
\cite{sid1},\cite{sid2}.

We start with a initial vector $x_0$ and the method we interested in
is to generate a sequence of vectors $x_1,x_2,\cdots,$ which
satisfies
$$x_m=x_0+q_{m-1}(A)r_0, \ \ r_0=:b-Ax_0,$$
where $q_{m-1}(\lambda)=\sum\limits_{i=1}^{m-a}c_i\lambda^{a+i-1}$
$(a:=ind(A))$. Then
\begin{equation}\label{2}
r_m:=b-Ax_m=p_m(A)r_0
\end{equation}
where $p_m(\lambda)=1-\lambda
q_{m-1}(\lambda)=1-\sum\limits_{i=1}^{m-a}c_i\lambda^{a+i}$.

The Krylov subspace we will use is
$$K_m(A;A^ar_0)=span\{A^ar_0,A^{a+1}r_0,\cdots,A^{m-1}r_0\}$$.
The vector $x_m$ produced by DGMRES satisfies
\begin{equation}\label{4}
\|A^ar_m\|=\min_{x\in x_0+K_m(A;A^ar_0)} \|A^a(b-Ax)\|.
\end{equation}

Now we give restarting DGMRES Algorithm DGMRES(m).
\begin{alg} DGMRES(m)\\
1. Choose initial guess $x_0$, and compute $r_0=b-Ax_0$ and $A^ar_0$;\\
2. Compute $\beta=\|A^ar_0\|$ and set $v_1=\beta^{-1}A^ar_0$;\\
3.Orthogonalize the vectors: $A^ar_0,A^{a+1}r_0,\cdots,A^mr_0$ by
Arnoldi-Gram-Schmidt process.\\

For $j=1,2,\cdots$ do

For $i=1,2,\cdots$ do

Compute $h_{i,j}=\langle AV_j,V_i\rangle$.

Compute $\widehat{v}_j=Av_j-\sum\limits_{i=1}^j h_{i,j}v_i$.

Set $h_{j+1,j}=\|\widehat{v}_i\|$ and
$v_{j+1}=\widehat{v}_j/h_{j+1,j}$.\\
4. For $k=1,2,\cdots,$, form the matrices $\widehat{V}_k\in
\Bbb{C}^{n\times k}$ and $\overline{H}_k\in \Bbb{C}^{(k+1)\times
k}$.
$$\widehat{V}_k=[v_1,v_2,\cdots,v_k];$$
$$\overline{H}_k=\left(\begin{array}{cccc}
h_{11} & h_{12} &\cdots &h_{1k}\\
h_{21} & h_{22} &\cdots &h_{2k}\\
0      & h_{32} &\cdots &h_{3k}\\
       &        &\cdots &      \\
\vdots & \ddots &\ddots &\vdots      \\
0      & \cdots &0      &h_{k+1,k}
\end{array}
\right)$$\\
5. For $m=a+1,\cdots$, form the matrix
$\widehat{H}_m=\overline{H}_m\overline{H}_{m-1}\cdots\overline{H}_{m-a}$.\\
6. Compute the QR factorization of $\widehat{H}_m$:
$\widehat{H}_m=Q_mR_m$, where $R_m$ is upper triangular.\\
7. Solve the system $R_my_m=\beta(Q_m^\ast e_1)$, where
$e_1=(1,0,\cdots,0)^T$.\\
8. Compute $x_m=x_0+\widehat{V}_{m-a}y_m$.\\
9. Restart: if $r_m=b-Ax_m$ satisfied with the residual norm then
stop, else let $x_0=x_m$ and go to 2.
\end{alg}

In \cite{sid1,sid2}, the convergence analysis of DGMRES is given as
follows.
\begin{lem}\cite{sid1,sid2}
Denote the spectrum of $A$ by $\sigma(A)$ and choose  $\Omega$ to be
a closed domain in the $\lambda-$plane that contains
$\sigma(A)\backslash\{0\}$ but not $\lambda=0$, such that its
boundary is twice differentiable with respect to arc-length. Denote
by $\Phi(\lambda)$ the conformal mapping of the exterior of $\Omega$
onto the exterior of the unit disk $\{w:|w|\geq 1\}$. Then the
vector $x_m$ extracted from $K_m(A;A^ar_0)$ satisfies
\begin{equation}\label{5}
\|A^ar_m\|\leq Km^{a+2(\widehat{k}-1)}\rho^m,
\end{equation}
for all $m$, where $K$ is a positive constant independent of $m$,
$\widehat{k}=max\{k_j:k_j=ind(A-\lambda_j),\lambda_j\in\sigma(A)\
\backslash \{0\}\}$ and $\rho=|1/\Phi(0)|<1$.
\end{lem}

We know the $x_m's$ produced by DGMRES is
$$x_m=x_0+\sum\limits_{i=1}^{m-a}c_iA^{a+i-1}r_0=x_0+p(A)r_0.$$
Let $z_1,z_2,\cdots,z_l$ is a set of linear independent eigenvectors
corresponding to eigenvalues $\lambda_1,\lambda_2,\cdots,\lambda_l$.
We add some new vectors
$\widehat{z}_{l+1},\widehat{z}_{l+2},\cdots,\widehat{z}_n$ such that
$z_1,\cdots z_l,\widehat{z}_{l+1},\cdots,\widehat{z}_n$ is a basis
in $\Bbb{C}^n$. Then
$r_0=b-Ax_0=\sum\limits_{i=1}^l\beta_iz_i+\sum\limits_{i=1+1}^n\beta_i\widehat{z}_i$
and
\begin{eqnarray*}
&&r_m=b-Ax_m=b-A(x_0+p(A)r_0)\\
&=&r_0-Ap(A)[\sum\limits_{i=1}^l\beta_iz_i+\sum\limits_{i=l+1}^n\beta_i\widehat{z}_i]\\
&=&r_0-\sum\limits_{i=1}^l\lambda_ip(\lambda_i)z_i-\sum\limits_{i=l+1}^n\beta_iAp(A)\widehat{z}_i\\
&=&\sum\limits_{i=1}^l\beta_iz_i+\sum\limits_{i=l+1}^n\beta_i\widehat{z}_i-\sum\limits_{i=1}^l
\lambda_ip(\lambda_i)z_i-\sum\limits_{i=l+1}^n\beta_iAp(A)\widehat{z}_i\\
&=&\sum\limits_{i=1}^l(\beta_iz_i-\lambda_ip(\lambda_i)z_i)+
\sum\limits_{i=l+1}^n(\beta_i-\beta_iAp(A))\widehat{z}_i.
\end{eqnarray*}
Thus
\begin{eqnarray*}
A^ar_m&=&\sum\limits_{i=1}^l[\beta_i\lambda_i^az_i-\lambda_i^{a+1}p(\lambda_i)z_i]
+\sum\limits_{i=l+1}^n[\beta_i\lambda_i^a\widehat{z}_i-\beta_ip(A)A^{a+1}\widehat{z}_i].
\end{eqnarray*}

From Lemma 2, we know
\begin{eqnarray}\label{6}
&&\|\sum\limits_{i=1}^l[\beta_i\lambda_i^az_i-\lambda_i^{a+1}p(\lambda_i)z_i]+
\sum\limits_{i=l+1}^n[\beta_i\lambda_i^a\widehat{z}_i-\beta_ip(A)A^{a+1}\widehat{z}_i]\|\\
\nonumber
 &\leq&Km^{a+2(\widehat{k}-1)}\rho^m.
\end{eqnarray}

\section{DGMRES augmented with eigenvectors}
In order to accelerate the convergence of restarting GMRES, Morgan
suggested some eigenvectors corresponding to a few of the smallest
eigenvalues add to the Krylov subspace for GMRES. Then the
convergence can be much faster. This method can be used in the
process of DGMRES.

Let $k$ be the number of the eigenvectors added to the subspace. Let
$r_0$ be the initial residual vector and $a=ind(A)$.
$K_m:=span\{A^ar_0,A^{a+1}r_0,\cdots,A^{m-1}r_0\}$ and
$K_{m,k}=span\{A^ar_0,A^{a+1}r_0,\cdots,A^{m-1}r_0,z_1,z_2,\cdots,z_k\}$,
where $z_1,z_2,\cdots,z_k$ are the eigenvectors added to $K_m$. The
approximated Drazin-inverse solution will be extracted from
$K_{m,k}$. That is
$$x_m=x_0+\sum\limits_{i=1}^k\alpha_iz_i+\sum\limits_{i=1}^{m-a}c_iA^{a+i-1}r_0$$.
which satisfies
\begin{equation}\label{11}
\|A^ar_m\|=\min_{x\in x_0+K_{m,k}(A;A^ar_0)} \|A^a(b-Ax)\|.
\end{equation}

Suppose $z_1,z_2,\cdots,z_k,\cdots,z_l$ is the set of linear
independent eigenvectors and
$z_1,\cdots,z_k,\cdots,z_l,\widehat{z}_{l+1},\cdots,\widehat{z}_n$
is a basis of $\Bbb{C}^n$. So we can write
$$r_0=\sum\limits_{i=1}^k\beta_iz_i+\sum\limits_{i=k+1}^l\beta_iz_i+\sum\limits_{i=l+1}^n\beta_i\widehat{z}_i,$$
and
\begin{eqnarray*}
&&x_m=x_0+\sum\limits_{i=1}^k\alpha_iz_i+p(A)r_0\\
&=&x_0+\sum\limits_{i=1}^k\alpha_iz_i+p(A)[\sum\limits_{i=1}^l\beta_iz_i+\sum\limits_{i=l+1}^n\beta_i\widehat{z}_i].
\end{eqnarray*}
Thus
\begin{eqnarray*}
&&r_m=b-Ax_m\\
&=&b-Ax_0-\sum\limits_{i=1}^k\alpha_i\lambda_iz_i-\sum\limits_{i=1}^l\lambda_ip(\lambda_i)\beta_iz_i-\sum\limits_{i=l+1}^nAp(A)\beta_i\widehat{z}_i\\
&=&r_0-\sum\limits_{i=1}^k\alpha_i\lambda_iz_i-\sum\limits_{i=1}^k\lambda_ip(\lambda_i)\beta_iz_i-\sum\limits_{i=k+1}^l\lambda_ip(\lambda_i)\beta_iz_i
-\sum\limits_{i=l+1}^n\beta_iAp(A)\widehat{z}_i
\end{eqnarray*}
and
\begin{eqnarray*}
&&A^ar_m=A^ar_0-\sum\limits_{i=1}^k\alpha_i\lambda_i^{a+1}z_i-\sum\limits_{i=1}^k\lambda_i^{a+1}p(\lambda_i)\beta_iz_i\\
&&\ \
-\sum\limits_{i=k+1}^l\lambda_i^{a+1}p(\lambda_i)\beta_iz_i-\sum\limits_{i=l+1}^n\beta_iAp(A)\widehat{z}_i\\
&=&\sum\limits_{i=1}^k\beta_i\lambda_i^az_i+\sum\limits_{i+1}^l\beta_i\lambda_i^az_i+\sum\limits_{i=l+1}^n\beta_iA^a\widehat{z}_i\\
&&\ \
-\sum\limits_{i=1}^k\alpha_i\lambda_i^{a+1}z_i-\sum\limits_{i=1}^k\lambda_i^{a+1}p(\lambda_i)\beta_iz_i\\
&&\ \
-\sum\limits_{i=k+1}^l\lambda_i^{a+1}p(\lambda_i)\beta_iz_i-\sum\limits_{i=l+1}^n\beta_iAp(A)\widehat{z}_i\\
&=&\sum\limits_{i=1}^k[\beta_i\lambda_i^a-\alpha_i\lambda_i^{a+1}-\lambda_i^{a+1}p(\lambda_i)\beta_i]z_i\\
&&+\sum\limits_{i=k+1}^l(\beta_i\lambda_i^a-\lambda_i^{a+1}p(\lambda_i)\beta_i)z_i+\sum\limits_{i=l+1}^n
[\beta_iA^a-\beta_iAp(A)]\widehat{z}_i\\
\end{eqnarray*}
Since DGMRES augmented with eigenvectors makes $\|A^ar_m\|$
minimized, we have
$\alpha_i=\frac{\beta_i\lambda_i^a-\lambda_i^{a+1}p(\lambda_i)\beta_i}{\lambda_i^{a+1}}$
and
\begin{equation}\label{7}
A^ar_m=\sum\limits_{i=k+1}^l(\beta_i\lambda_i^a-\lambda_i^{a+1}p(\lambda_i)\beta_i)z_i
+\sum\limits_{i=l+1}^n[\beta_iA^a-\beta_iAp(A)]\widehat{z}_i
\end{equation}

From Lemma 1, (\ref{5}), and (\ref{7}), we get the follow theorem.
\begin{thm}
Denote the spectrum of $A$ by $\sigma
(A)=\{\lambda_1,\lambda_2,\cdots,\lambda_l\}$. $\sigma_a
(A):=\{\lambda_{k+1},\lambda_{k+2},\cdots,\lambda_l\}$ $(k\leq l)$.
 Choose  $\Omega_a$ to be
a closed domain in the $\lambda-$plan that contains
$\sigma_a(A)\backslash\{0\}$ but not $\lambda=0$, such that its
boundary is twice differentiable with respect to arc-length. Denote
by $\Phi_a(\lambda)$ the conformal mapping of the exterior of
$\Omega_a$ onto the exterior of the unit disk $\{w:|w|\geq 1\}$.
Then the vector $x_m$ generated by $DGMRES$ augmented with
eigenvectors satisfies
\begin{equation}\label{8}
\|A^ar_m\|\leq Km^{a+2(\widehat{k}-1)}\rho_a^m,
\end{equation}
for all $m$, where $K$ is a positive constant independent of $m$,
$\widehat{k}=max\{k_j:k_j=ind(A-\lambda_j),\lambda_j\in\sigma_a(A)\
\backslash \{0\}\}$ and $\rho_a=|1/\Phi_a(0)|<1$.
\end{thm}

Compared (\ref{8}) with (\ref{5}), we know the convergence of DGMRES
augmented with eigenvectors is faster than that of DGMRES.

The implementation of DGMRES augmented with eigenvectors is
different from that of GMRES. We derive the algorithm step by step.

1. Let the initial vector $x_0=0$. Compute $\beta=\|A^ar_0\|$ and
set $v_1=(A^ar_0)/\beta$.

2. Orthogonalize the vectors $A^ar_0,\cdots,A^{m-1}r_0$ via
Arnoldi-Gram-Schmidt process.

 For $j=1,2,...,m-a$.

 For $i=1,2,\cdots,j$

 Compute $h_{ij}=\langle Av_j,v_i\rangle$.

 $\widehat{v_{j+1}}=Av_j-h_{ij}v_i$

 end i

 Let $h_{j+1,j}=\|\widehat{v}_j\|$

 Set $v_{j+1}=\widehat{v}_j/h_{j+1,j}$

 end j\\

 Consequently, we get a orthogonal matrix $V_{m-a}=[v_1,v_2,\cdots,
 v_{m-a}]$ and a $(m-a+1)\times (m-a)$ Hessenberg matrix
 $\overline{H}^{\sharp}$.

 3. Compute k approximated smallest  magnitude nonzero  eigenvalues of $A$ and add their corresponding eigenvectors $z_1,\cdots,z_k$ to the
 subspace. Let $H^{\sharp}=[h_{ij}]_{(m-a)\times (m-a)}$. Obviously
 $V_{m-a}^TAV_{m-a}=H^{\sharp}$. When $m\leq q-1$, where $q$ is the
 degree of the minimal polynomial of $A$, the $\overline{H}^{\sharp}$
 has full rank. We suppose $\lambda_1,\lambda_2,\cdots,\lambda_k$
 are k smallest eigenvalues of $H^{\sharp}$ and $y_1,y_2,\cdots,y_k$
 are the eigenvectors corresponding to them. Then $z_i:=V_{m-a}y_i$
 $(i=1,2,\cdots,k)$ are the approximated eigenvectors of $A$. We add
 them to the subspace and orthogonalize.

For $j=m-a+1,\cdots,m-a+k$

For $i=1,2,\cdots,j$

Compute $h_{ij}=\langle Az_{j-m+a},v_i\rangle$.

Compute $\widehat{v}_{j+1}=Az_{j-m+a}-h_{ij}v_i$

end i

Set $h_{j+1,j}=\|\widehat{v}_{j+1}$

Set $v_{j+1}=\widehat{v}_{j+1}/h_{j+1,j}$

end j.

Denote $W=[v_1,v_2,\cdots,v_{m-a},z_1,\cdots,z_k]$,
$\overline{H}^{(0)}=[h_{ij}]_{m-a+k+1,m-a+k}$, and
$V^{(0)}=[v_1,v_2,\cdots,v_{m-a},v_{m-a+1},\cdots,v_{m-a+k+1}]$. It
is easy to see $AW=V^{(0)}\overline{H}^{(0)}$ and $x_m=x_0+Wy_m$ for
some $y_m\in \Bbb{C}^{n\times n}$. It follows that
$r_m=r_0-AWy_m=r_0-V^{(0)}\overline{H}^{(0)}y_m$. So we get
\begin{equation}\label{8}
A^ar_m=A^ar_0-A^{a+1}Wy_m=\beta v_1-A^{a+1}Wy_m.
\end{equation}
From this, we should continue to deal with $A^{a+1}W$.
\begin{eqnarray*}\label{9}
A^{a+1}W&=&A^a(AW)=A^aV^{(0)}\overline{H}^{(0)} \\
&=&A^{a-1}(AV^{(0)})\overline{H}^{(0)} \\
&=&A^{a-1}[v_1,v_2,\cdots,v_{m-a},v_{m-a+1},v_{m-a+2},\cdots,v_{m-a+k+2}]\overline{H}^{(1)}\overline{H}^{(0)}\\
&=&A^{a-1}V^{(1)}\overline{H}^{(1)}\overline{H}^{(0)}\\
&\vdots&\\
&=&A^{a-(a-1)}V^{(a-1)}\overline{H}^{(a-1)}\cdots
\overline{H}^{(0)}\\
&=&V^{(a)}\overline{H}^{(a)}\overline{H}^{(a-1)}\cdots
\overline{H}^{(0)}\\
\end{eqnarray*}
Denote $\overline{H}^{(a)}\cdots \overline{H}^{(0)}$ by
$\overline{H}$. From the above, we get
\begin{equation}\label{10}
A^{a+1}W=V^{(a)}\overline{H}.
\end{equation}
Since $r_m$ is the residual vector produced by DGMRES augmented with
eigenvectros, from (\ref{8}) and (\ref{10}), we have
\begin{eqnarray*}
\|A^ar_m\|&=&\beta v_1-A^{a+1}Wy_m\|=\|\beta
v_1-V^{(a)}\overline(H)y_m\| \\
&=&\|\beta e_1-\overline{H}y_m\|=\min\limits_{y}\|\beta
e_1-\overline{H}y\|
\end{eqnarray*}
Apply the QR factorization of $\overline{H}$: $\overline{H}=QR$,
where $R$ is upper triangular. Thus $y_m$ satisfies
$Ry_m=\beta(Q^Te_1)$ and $x_m$ follows.

4. The practical computation of $\overline{H}$. From the above we
know
$$A[v_1,\cdots,v_{m-a},z_1,\cdots,z_k]=[v_1,\cdots,v_{m-a},v_{m-a+1},\cdots,v_{m-a+k+1}]\overline{H}^{(0)},$$
and
\begin{eqnarray*}
&&A[v_1,v_2,\cdots,v_{m-a},v_{m-a+1},\cdots,v_{m-a+k+1}\\
&=&[v_1,v-2,\cdots,v_{m-a},v_{m-a+1},v_{m-a+2}^{(1)},\cdot,v_{m-a+k+1}^{(1)}]\overline{H}^{(1)}
\end{eqnarray*}

where $$  \overline{H}^{(1)}=\left(
\begin{array}{ccccccc}
h_{11} &  h_{12}  &    \cdots   &   h_{1,m-a} &
\widehat{h}_{1,m-a+1}   &    \cdots    &
 \widehat{h}_{1,m-a+k+1}\\
h_{21}  &  h_{22}   &  \cdots   &h_{2,m-a}  &\widehat{h}_{2,m-a+1}
&\cdots  &\widehat{h}_{1,m-a+k+1}\\
0   &h_{32}   &\cdots   &h_{3,m-a}  &\widehat{h}_{3,m-a+1}   &\cdots
&\widehat{h}_{3,m-a+k+1}\\
      &    \cdots   &\cdots  & & & & \\
0 & 0 &\cdots &h_{m-a+1,m-a} &\widehat{h}_{m-a+1,m-a+1} &\cdots
&\widehat{h}_{m-a+1,m-a+k+1}\\
0 & 0 &\cdots &0 &\widehat{h}_{m-a+2,m-a+1}  &\cdots
&\widehat{h}_{m-a+2,m-a+k+1}\\
     &  \cdots   &\cdots &&&& \\
0 & 0 &\cdots &0  &0  &\cdots   &\widehat{h}_{m-a+k+2,m-a+k+1}
\end{array}
\right) $$

$\widehat{h}_{ij}$ is the new entries we need to compute when from
$\overline{H}^{(0)}$ to $\overline{H}^{(1)}$, that is, we need to
compute $\frac{[2(m-a+2)+k]\times (k+1)}{2}$ entries additional.

Next
\begin{eqnarray*}
&&A[v_1,v_2,\cdots, v_{m-a}, v_{m-a+1},
v^{(1)}_{m-a+2},\cdots,v^{(1)}_{m-a+k+1}]\\
&=&[v_1,\cdots,v_{m-a+1},v_{m-a+2}^{(1)},v_{m-a+3}^{(2)},\cdots,v_{m-a+k+2}^{(2)}]\overline{H}^{(2)}.
\end{eqnarray*}

where $$  \overline{H}^{(1)}=\left(
\begin{array}{ccccccc}
h_{11} &h_{12} &\cdots   &h_{1,m-a+1}  &\widehat{h}_{1,m-a+2}
&\cdots   &\widehat{h}_{1,m-a+k+2}\\
h_{21}   &h_{22}   &\cdots   &h_{2,m-a+1}  &\widehat{h}_{2,m-a+2}
&\cdots    &\widehat{h}_{2,m-a+k+2}\\
0 & h_{32}  &\cdots   &h_{3,m-a+1}    &\widehat{h}_{3,m-a+2} &\cdots
&\widehat{h}_{3,m-a+k+2}\\
& &\vdots &\vdots &&&\\
0 & 0 &\cdots  &h_{m-a+2,m-a+1}  &\widehat{h}_{m-a+2,m-a+2}  &\cdots
&\widehat{h}_{m-a+2,m-a+k+2}\\
& & \vdots &\vdots &&&\\
0 &0 &\cdots &0&0& \cdots& \widehat{h}_{m-a+k+2,m-a+k+2}\\
\end{array}
\right) $$

Similarly, from $\overline{H}^{(1)}$ to $\overline{H}^{(2)}$, there
are $\frac{[2(m-a)+5+k]\times (k+1)}{2}$ need to be computed.
Following this way, we can continue to get
$\overline{H}^{(3)},\cdots, \overline{H}^{(a)}$ and  $\overline{H}$
can be computed.

We summarize the above as the following algorithm.

\begin{alg} DGMRES augmented with eigenvectors\\
1. Pick initial vector $x_0$ and compute $r_0=b-Ax_0$ and $\beta=\|A^ar_0\|$.\\
2. Apply Arnoldi process to $A^ar_0,\cdots, A^{m-1}r_0$.

 $v_1=r_0/\|r_0\|$

 For $j=1,\cdots,m-a$

 For $i=1,\cdots,j$

 Compute $h_{ij}=\langle Av_j,v_i\rangle$

 Compute $\widehat{v}_{j+1}=Av_j-h_{ij}v_i$

 end i

 Set $h_{j+1,j}=\|v_j\|$

 $v_{j+1}=\widehat{v}_{j+1}/h_{j+1,j}$

 end j\\
 4. Denote $[h_{ij}]_{(m-a)\times (m-a)}$ by $H$ and compute its k
 eigenvectors $y_1,\cdots,y_k$ corresponding to k smallest magnitude eigenvalues.
 Compute $z_i=Vy_i$. We add $z_i's$ to the subspace and denote
 $W=[v_1,\cdots,v_{m-a},z_1,\cdots,z_k]$. We apply Arnoldi process
 to $AW$ and get the matrix $V^{(0)}$ and $\overline{H}^{(0)}$.

 For $j=m-a+1,\cdots,m-a+k$

 For $i=1,2\cdots,j$

 Compute $h_{ij}=\langle Az_{j-(m-a)},v_i\rangle$

 Compute $\widehat{v}_{j+1}=Az_{j-(m-a)}-h_{ij}v_i$

 end i

 Set $h_{j+1,j}=\|v_j\|$

 Set $v_{j+1}=\widehat{v}_{j+1}/h_{j+1,j}$

 end j\\
5. Compute $\overline{H}^{(1)},\cdots,\overline{H}^{(a)}$.

For $t=1,2,\cdots,a$

For $j=m-a+t,\cdots,m-a+k+t$

For $i=1,2,\cdots,j$

Compute $h_{ij}=\langle Av_j,v_i\rangle$

$\widehat{v}_{j+1}=Av_j-h_{ij}v_i$

end i

$h_{j+1,j}=\|v_j\|$

$v_{j+1}=\widehat{v}_{j+1}/h_{j+1,j}$

end j

Compute $\overline{H}=\overline{H}^{(a)}\cdots
\overline{H}^{(0)}$.\\
6. Compute the QR factorization of $\overline{H}$:
$\overline{H}=QR$, where $R$ is upper triangular. Set $c=\beta Q^T
e_1$. Solve the least-square problem $R_my_m=\beta (Q^Te_1)$ and
$x_m=x_0+Wy_m$. If $r_m:=b-Ax_m$ satisfies the residual norm then
stop; else  set $x_0=x_m$ and go to 2.
\end{alg}

\section{Numerical examples}
For convenience, We denote the DGMRES augmented with eigenvectors by
ADGMRES. In this section we present some numerical examples to
compare ADGMRES to DGMRES. All the experiments were performed in
$MATLAB^\circledR$ 7.5 on an Inter Core 2 Duo 2000MHz PC with main
memory 1000M.

{\bf Example 1} We take $A$ to be a 12 by 12 singular matrix which
has the following Jordan canonical form
$$\left(
\begin{array}{cccccccccccc}
1 &1 &0 &0 &0 &0 &0 &0 &0 &0 &0 &0\\
0 &1 &1 &0 &0 &0 &0 &0 &0 &0 &0 &0\\
0 &0 &1 &0 &0 &0 &0 &0 &0 &0 &0 &0\\
0 &0 &0 &3 &1 &0 &0 &0 &0 &0 &0 &0\\
0 &0 &0 &0 &3 &1 &0 &0 &0 &0 &0 &0\\
0 &0 &0 &0 &0 &3 &0 &0 &0 &0 &0 &0\\
0 &0 &0 &0 &0 &0 &7 &0 &0 &0 &0 &0\\
0 &0 &0 &0 &0 &0 &0 &8 &0 &0 &0 &0\\
0 &0 &0 &0 &0 &0 &0 &0 &9 &1 &0 &0\\
0 &0 &0 &0 &0 &0 &0 &0 &0 &9 &0 &0\\
0 &0 &0 &0 &0 &0 &0 &0 &0 &0 &0 &1\\
0 &0 &0 &0 &0 &0 &0 &0 &0 &0 &0 &0
\end{array}
\right)
$$
The right side $b=(1,1,\cdots,1)^T$. The ADGMRES uses $m=6,k=1$ and
DGMRES uses $m=7$ and denote them by DGMRES(6,1), DGMRES(7)
respectively. So they use the same size subspace. When
$\frac{\|A^ar_m\|}{\|A^ab\|}<\epsilon$, we stop. The convergence
curve for ADGMRES and DGMRES are indicated in Fig 1. It is easy to
see ADGMRES convergence faster than DGMRES

\begin{center}
\scalebox{0.6}{\includegraphics{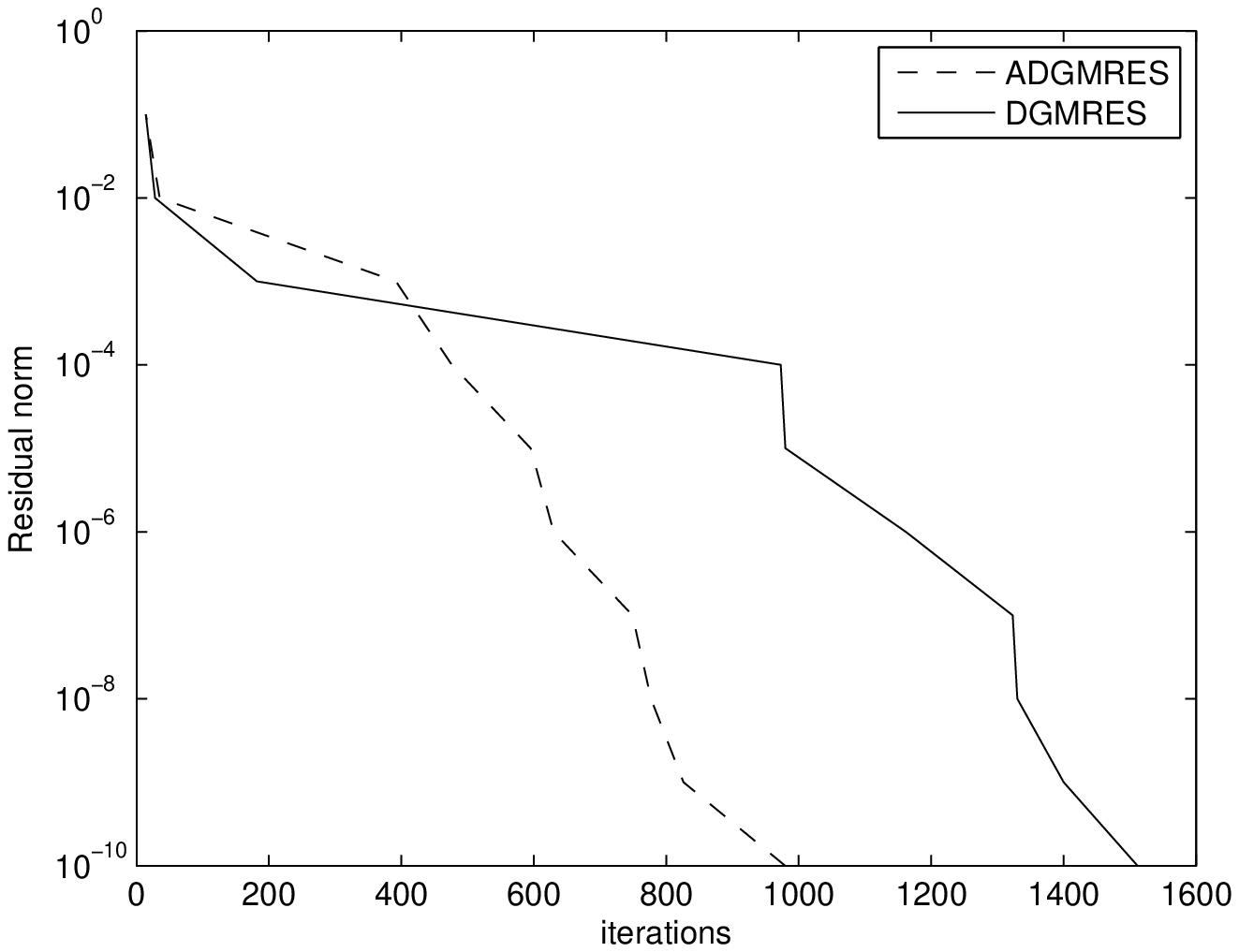}}

{\small Fig1. The convergence curves for ADGMRES(6,1) and DGMRES(7)}
\end{center}

{\bf Example 2.} In this example the matrix $A$ is the same as the
above but let $a_{7,7}=1000$. Obviously such $A$ has a larger
eigenvalue 1000.  We expect ADGMRES(4,1) convergence much faster
than $DGMRES(5)$. As indicated in Fig 2., it is the case.

\begin{center}
\scalebox{0.6}{\includegraphics{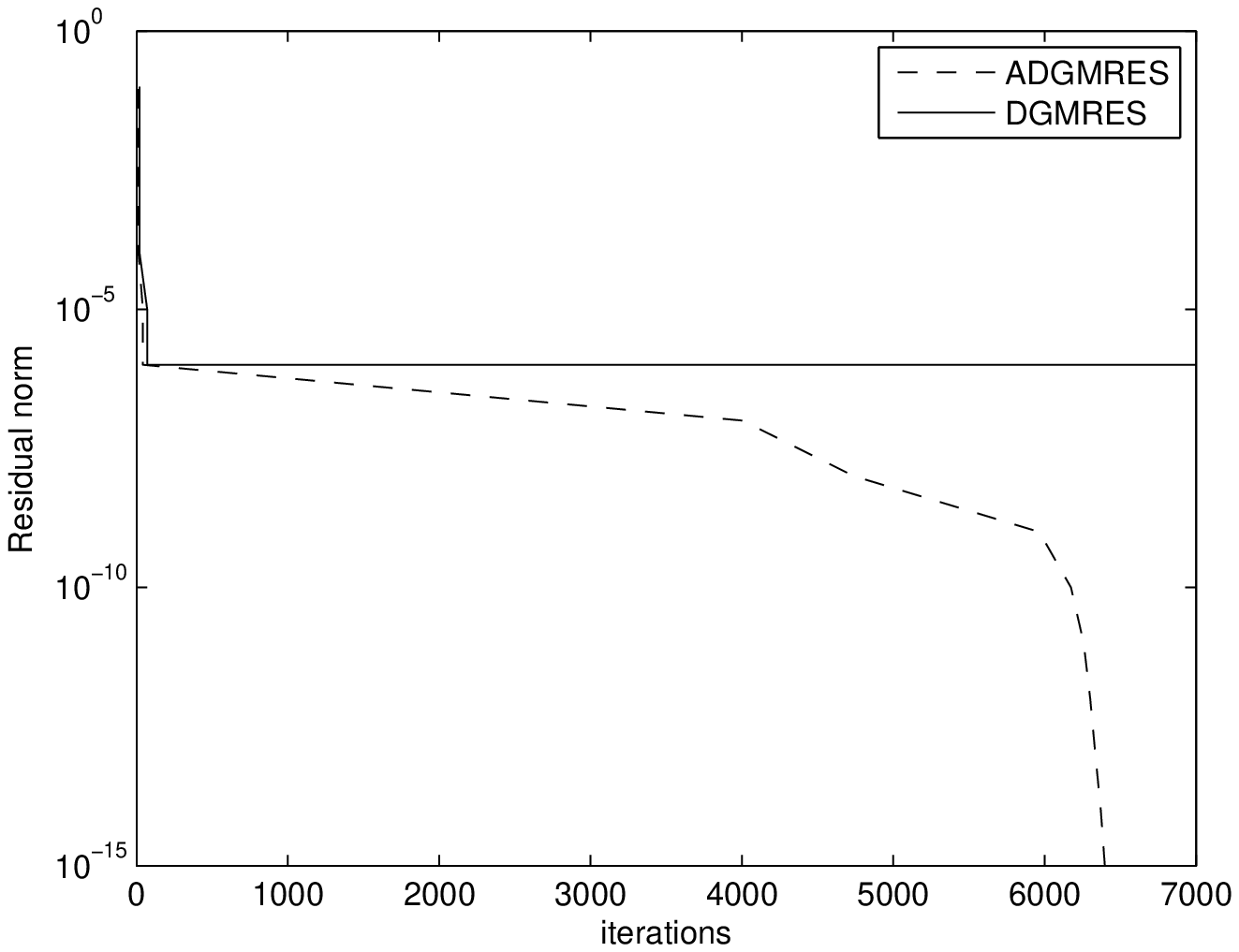}}

{\small Fig2. The convergence curves for ADGMRES(4,1) and DGMRES(5)}
when $A$ has a very lager eigenvalue.
\end{center}

From Fig 2. we can see after 70 iterations the DGMRES(5) stagnates
but after 5000 iterations the curve of ADGMRES(4,1) decreases
steeply.

{\bf Example 3.} In this example we take $a_{77}=0.001$, that is,
$A$ has a very small eigenvalue. From Fig 3. we see DGMRES(5)
stagnates after 50 iterations but ADGMRES(4,1) still works well.

\begin{center}
\scalebox{0.6}{\includegraphics{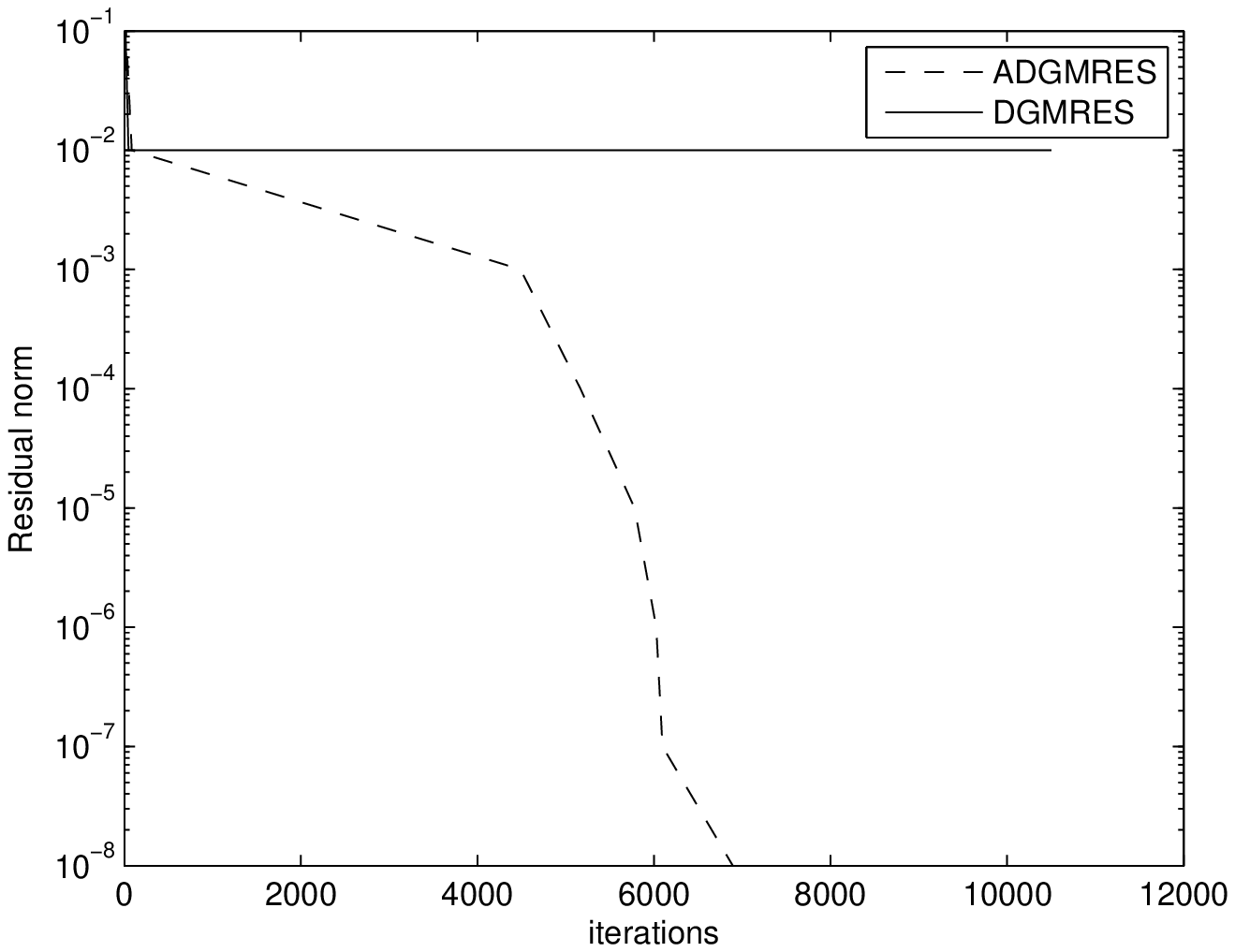}}

{\small Fig3. The convergence curves for ADGMRES(4,1) and DGMRES(5)}
when $A$ has a very small eigenvalue.
\end{center}

{\bf Example 4.} This example comes from \cite{zhou2}.
$$A=\left(
\begin{array}{cccc}
1 & 1 &1 &2\\
0 &1 &3 &4\\
0 &0 &1 &1\\
0 &0 &0 &0
\end{array}
\right), b=\left(\begin{array}{c}
 -4\\
7\\
1\\
0\end{array}\right)$$

In \cite{zhou2} The author found that DGMRES(2) converges faster
than DGMRES(3). Since $A$ has only one nonzero eigenvalue $1$, we
infer from Theorem that the convergence of ADGMRES(2,1) is as the
same as that of DGMRES(2). From  Table 1 we observe it is the case.
At the same time we also see DGMRES(3) stagnates. The residual norm
produced by the three methods is indicated in the following table.

\begin{center}
\begin{tabular}{c| c c c}
runs     &ADGMRES(2,1)         &DGMRES(3)        &DGMRES(2)\\
\hline

50 &0.0038      & 0.00279         &0.0038\\
100 &$1.23\times 10^{-5}$  &0.00276   &$1.23\times 10^{-5}$\\
200  &$1.71\times 10^{-9}$  &0.00276   &$1.72\times 10^{-9}$\\
300  &$6.155\times 10^{-14}$   &0.00276  &$6.145\times 10^{-14}$\\
\end{tabular}

{\small Table 1. The residual norm produced by ADGMRES(2,1),
DGMRES(3), DGMRES(2) after 50, 100, 200, 300 runs respectively.}
\end{center}

\section{Conclusion.} The DGMRES method augmented with eigenvectors
can improve the convergence, especially when the matrix has small or
large eigenvalues. The DGMRES method often stagnates (see
\cite{zhou1,zhou2}). The DGMRES augmented with eigenvectors is a
good choice to overcome the stagnation. When $k$, the number of the
eigenvectors added to the subspace, is large, the method is
expensive. So in practical use we usually choose a small $k$.

\bibliographystyle{amsplain}

\end{document}